\documentclass{article}
\usepackage{amssymb}
\usepackage{amsmath}

\newtheorem{theorem}{Theorem}

\newtheorem{corollary}[theorem]{Corollary}

\newtheorem{definition}[theorem]{Definition}

\newtheorem{proposition}[theorem]{Proposition}
\newtheorem{remark}[theorem]{Remark}

\newenvironment{proof}[1][Proof]{\noindent\textbf{#1.} }{\ \rule{0.5em}{0.5em}}

\newcommand{\RR}{\mathbb{R}}

\renewcommand{\a}{\alpha}
\renewcommand{\b}{\beta}

\begin{document}

\title{On Opial-type inequalities via fractional calculus}

\author{\textbf{Portilla, Ana}\\Saint Louis University, Madrid Campus\\ Avenida del Valle 34, 28003 Madrid, Spain\\\texttt{aportil2@slu.edu}\\
\and \textbf{Rodr\'{\i}guez, Jos\'{e} M.}\\Universidad Carlos III de Madrid, Departamento de Matem\'aticas\\Avenida de la Universidad 30, 28911 Legan\'{e}s, Madrid, Spain\\\texttt{jomaro@math.uc3m.es}\\
\and \textbf{Sigarreta Almira, Jos\'{e} M.}\\Universidad Aut\'{o}noma de Guerrero, Centro Acapulco\\CP 39610, Acapulco de Ju\'{a}rez, Guerrero, Mexico\\\texttt{josemariasigarretaalmira@hotmail.com} }

\maketitle

\begin{abstract}
Inequalities play an important role in pure and applied mathematics.
In particular, Opial inequality
plays a main role in the study of the existence and uniqueness of initial and boundary value
problems for differential equations.
It has several interesting generalizations.
In this work we prove some new Opial-type inequalities, and we apply them to generalized Riemann-Liouville-type
integral operators.
\end{abstract}

\bigskip



\bigskip

\section{Introduction}

Integral inequalities are used in countless
mathematical problems such as approximation theory and spectral analysis, statistical analysis and the
theory of distributions. Studies involving integral inequalities play an important role in several areas
of science and engineering.

In recent years there has been a growing interest in the study of many classical inequalities applied to integral operators associated with different types of fractional derivatives, since integral inequalities and their applications
play a vital role in the theory of differential equations and applied mathematics.
Some of the inequalities studied are Gronwall, Chebyshev, Jensen-type, Hermite-Hadamard-type, Ostrowski-type, Gr\"uss-type, Hardy-type,
Gagliardo-Nirenberg-type, reverse Minkowski and reverse H\"older
inequalities (see, e.g., \cite{Dahmani,Han,Mubeen,Nisar,QRS,Rahman,Rahman2,Rashid,Sawano,Set}).

In this work we obtain new Opial-type inequalities, and we apply them to the generalized Riemann-Liouville-type
integral operators defined in \cite{BCRS}, which include most of known Riemann-Liouville-type
integral operators.

\section{Preliminaries}

One of the first operators that can be called fractional is the Riemann-Liouville fractional derivative of order $\alpha \in \mathbb{C}$, with $Re(\alpha)> 0$, defined as follows (see \cite{GM}).

\begin{definition} \label{d:RL}
	Let $a < b$ and $f \in L^{1}((a,b);\mathbb{R})$.
	The \emph{right and left side Riemann-Liouville fractional integrals of order} $\alpha$, with $Re(\alpha)> 0$, are defined, respectively, by
	\begin{equation}\label{e:RL+}
		_{  }^{ RL }\!{ { J }_{ { a }^{ + } }^{ \alpha  } }f(t)=\frac { 1 }{ \Gamma (\alpha ) } \int _{ a }^{ t }{ { (t-s) }^{ \alpha -1 }f(s)\,ds },
	\end{equation}
	and
	\begin{equation}\label{e:RL-}
		_{  }^{ RL }\!{ { J }_{ { b }^{ - } }^{ \alpha  } }f(t)=\frac { 1 }{ \Gamma (\alpha ) } \int _{ t }^{ b }{ { (s-t) }^{ \alpha -1 }f(s)\,ds },
	\end{equation}
	with $t \in (a,b)$.
\end{definition}

When $\a \in (0,1)$, their corresponding \emph{Riemann-Liouville fractional derivatives} are given by
$$
\begin{aligned}
	\big({ _{  }^{ RL }\!D_{ a^{ + } }^{ \alpha  } }f\big)(t)
	& =\frac { d }{ dt } \left( { _{  }^{ RL }\!{ { J }_{ { a }^{ + } }^{ 1-\alpha  } }f(t) } \right)
	=\frac { 1 }{ \Gamma (1-\alpha ) } \, \frac { d }{ dt } \int _{ a }^{ t }{ { \frac { f(s) }{ (t-s)^{ \alpha  } }  }\,ds },
	\\
	\big({ _{  }^{ RL }\!D_{ b^{ - } }^{ \alpha  } }f\big)(t)
	& =-\frac { d }{ dt } \left( { _{  }^{ RL }\!{ { J }_{ { b }^{ - } }^{ 1-\alpha  } }f(t) } \right)
	=-\frac { 1 }{ \Gamma (1-\alpha ) } \,\frac { d }{ dt } \int _{ t }^{ b }{ { \frac { f(s) }{ (s-t)^{ \alpha  } }  }\,ds }.
\end{aligned}
$$

Other definitions of fractional operators are the following ones.

\begin{definition} \label{d:H}
	Let $a < b$ and $f \in L^{1}((a,b);\mathbb{R})$.
	The \emph{right and left side Hadamard fractional integrals of order} $\alpha$, with Re($\alpha)>0$, are defined, respectively, by
	\begin{equation}\label{e:H+}
		{ { H }_{ { a }^{ + } }^{ \alpha  } }f(t)=\frac { 1 }{ \Gamma (\alpha ) } \int _{ a }^{ t }{ { \Big(\log\frac { t }{ s } \,\Big) }^{ \alpha -1 }\frac { f(s) }{s } \, ds } ,
	\end{equation}
	and
	\begin{equation}\label{e:H-}
		{ { H }_{ { b }^{ - } }^{ \alpha  } }f(t)=\frac { 1 }{ \Gamma (\alpha ) } \int _{ t }^{ b }{ { \Big(\log\frac { s }{t } \,\Big) }^{ \alpha -1 }\frac { f(s) }{ s } \, ds } ,
	\end{equation}
	with $t \in (a,b)$.
\end{definition}

When $\a \in (0,1)$, \emph{Hadamard fractional derivatives} are given by the following expressions:
$$
\begin{aligned}
	{ \big( {  }^{ H }\!{ D_{ a^{ + } }^{ \alpha  } } }f\big)(t)
	& =t\,\frac { d }{ dt } \big( { { H }_{ { a }^{ + } }^{ 1-\alpha  } }f(t) \big)
	=\frac { 1 }{ \Gamma(1-\alpha ) } \, t \, \frac { d }{ dt }\int _{ a }^{ t }{ { \Big(\log \frac { t }{ s } \Big) }^{ -\alpha }\frac { f(s) }{ s } \,ds } ,
	\\
	{ \big( {  }^{ H }\!{ D_{ b^{ - } }^{ \alpha  } } }f\big)(t)
	& =-t\,\frac { d }{ dt } \big( { { { H }_{ { b }^{ - } }^{ 1-\alpha  } }f(t) } \big)
	= \frac {-1}{ \Gamma (1-\alpha ) } \, t \, \frac { d }{ dt } \int _{ t }^{ b }{ { \Big(\log\frac { s }{ t } \Big) }^{ -\alpha }\frac { f(s) }{ s } \,ds } ,
\end{aligned}
$$
with $t \in (a,b)$.

\medskip

\begin{definition} \label{d:20}
	Let $0 < a < b$,
	$g:[a, b] \rightarrow \mathbb{R}$ an increasing positive function on $(a,b]$ with continuous derivative on $(a,b)$,
	$f:[a, b]\rightarrow \mathbb{R}$ an integrable function, and $\alpha \in (0,1)$ a fixed real number.
	The right and left side fractional integrals in \cite{KMS} of order $\alpha$ of $f$ with respect to $g$ are defined,
	respectively, by
	\begin{equation}\label {e:fg+}
		I_{ g,a^+ }^{ \alpha  }f(t)=\frac { 1 }{ \Gamma (\alpha ) } \int _{ a }^{ t } \frac { g'(s)f(s) }{ { \big( g(t)-g(s) \big)  }^{ 1-\alpha  } } \,ds,
	\end{equation}
	and
	\begin{equation}\label {e:fg-}
		I_{g,b^- }^{ \alpha  }f(t)=\frac { 1 }{ \Gamma (\alpha ) } \int _{ t }^{ b } \frac { g'(s)f(s) }{ { \big( g(s)-g(t) \big)  }^{ 1-\alpha  } } \,ds,
	\end{equation}
	with $t \in (a,b)$.
\end{definition}

There are other definitions of integral operators in the global case, but they are slight modifications of the previous ones.

\section{General fractional integral of Riemann-Liouville type}

Now, we give the definition of a general fractional integral in \cite{BCRS}.

\begin{definition} \label{d:01}
	Let $a<b$ and $\a \in \RR^+$.
	Let $g:[a, b] \rightarrow \mathbb{R}$ be a positive function on $(a,b]$ with continuous positive derivative on $(a,b)$,
	and $G:[0, g(b)-g(a)]\times (0,\infty)\rightarrow \mathbb{R}$ a continuous function which is positive on $(0, g(b)-g(a)]\times (0,\infty)$.
	Let us define the function $T:[a, b] \times [a, b] \times (0,\infty)\rightarrow \mathbb{R}$ by
	$$
	T(t,s,\alpha)
	= \frac{ G\big( |g(t)-g(s)| ,\alpha \big) }{ g'(s) } \,.
	$$
	The \emph{right and left integral operators}, denoted respectively by $J_{T,a^+}^\a$ and $J_{T,b^-}^\a$,
	are defined for each measurable function $f$ on $[a,b]$ as
	\begin{equation}\label {e:oig+}
		J_{T,a^+}^\a f(t)=\int _{ a }^{ t } \frac { f(s) }{ T(t,s,\alpha ) } \, ds,
	\end{equation}
	\begin{equation}\label {e:oig-}
		J_{T,b^-}^\a f(t)=\int _{ t }^{ b } \frac { f(s) }{ T(t,s,\alpha ) } \, ds,
	\end{equation}
	with $t \in [a,b]$.
	
	We say that $f \in L_T^1[a,b]$ if $J_{T,a^+}^\a |f|(t), J_{T,b^-}^\a |f|(t) < \infty$ for every $t \in [a,b]$.
\end{definition}

\

Note that these operators generalize the integral operators in Definitions \ref{d:RL}, \ref{d:H} 
and \ref{d:20}:

\medskip

$(A)$ If we choose
$$
g(t)=t,
\quad G( x ,\alpha ) = \Gamma(\a)\, x^{ 1-\alpha},
\quad T(t,s,\alpha )=\Gamma (\alpha )\, |t-s|^{ 1-\alpha},
$$
then $J_{T,a^+}^\a$ and $J_{T,b^-}^\a$
are the right and left Riemann-Liouville fractional integrals $_{  }^{ RL }\!{ { J }_{ { a }^{ + } }^{ \alpha  } }$ and $_{  }^{ RL }\!{ { J }_{ { b }^{- } }^{ \alpha  } }$
in \eqref{e:RL+} and \eqref{e:RL-}, respectively.
Its corresponding right and left Riemann-Liouville fractional derivatives are
$$
\begin{aligned}
	\big({ _{  }^{ RL }\!D_{ a^{ + } }^{ \alpha  } }f\big)(t)
	=\frac { d }{ dt } \left( { _{  }^{ RL }\!{ { J }_{ { a }^{ + } }^{ 1-\alpha  } }f(t) } \right),
	\quad
	\big({ _{  }^{ RL }\!D_{ b^{ - } }^{ \alpha  } }f\big)(t)
	=-\frac { d }{ dt } \left( { _{  }^{ RL }\!{ { J }_{ { b }^{ - } }^{ 1-\alpha  } }f(t) } \right).
\end{aligned}
$$

\medskip

$(B)$ If we choose
$$
g(t)= \log t,
\quad G( x ,\alpha ) = \Gamma(\a)\, x^{ 1-\alpha},
\quad T(t,s,\alpha )=\Gamma (\alpha )\, t \, \Big|\log \frac{t}{s}\Big|^{ 1-\alpha},
$$
then $J_{T,a^+}^\a$ and $J_{T,b^-}^\a$
are the right and left Hadamard fractional integrals ${ H }_{ { a }^{ + } }^{\alpha}$ and ${ H }_{ { b }^{ - } }^{\alpha}$
in \eqref{e:H+} and \eqref{e:H-}, respectively.
Its corresponding right and left Hadamard fractional derivatives are
$$
\begin{aligned}
	{ \big( {  }^{ H }\!{ D_{ a^{ + } }^{ \alpha  } } }f\big)(t)
	=t\,\frac { d }{ dt } \big( { { H }_{ { a }^{ + } }^{ 1-\alpha} }f(t) \big) ,
	\quad
	{ \big( {  }^{ H }\!{ D_{ b^{ - } }^{ \alpha  } } }f\big)(t)
	=-t\,\frac { d }{ dt } \big( { { { H }_{ { b }^{ - } }^{ 1-\alpha} }f(t) } \big).
\end{aligned}
$$

\medskip

$(C)$ If we choose a function $g$ with the properties in Definition \ref{d:01} and
$$
G( x ,\alpha ) = \Gamma(\a)\, x^{ 1-\alpha},
\quad T(t,s,\alpha ) = \Gamma(\a)\, \frac{\, |\, g(t)-g(s) |^{ 1-\alpha}}{g'(s)} \,,
$$
then $J_{T,a^+}^\a$ and $J_{T,b^-}^\a$
are the right and left 
fractional integrals $I_{ g,a^+ }^{ \alpha  }$ and $I_{ g,b^- }^{ \alpha  }$
in \eqref{e:fg+} and \eqref{e:fg-}, respectively.

\medskip

\begin{definition} \label{d:0}
	Let $a<b$ and $\a \in \RR^+$.
	Let $g:[a, b] \rightarrow \mathbb{R}$ be a positive function on $(a,b]$ with continuous positive derivative on $(a,b)$,
	and $G:[0, g(b)-g(a)]\times (0,\infty)\rightarrow \mathbb{R}$ a continuous function which is positive on $(0, g(b)-g(a)]\times (0,\infty)$.
	For each function $f \in L_T^1[a,b]$, its \emph{right and left generalized derivative of order} $\alpha $ are defined, respectively,  by
	\begin{equation} \label{e:d0}
		\begin{aligned}
			D_{T,a^+ }^{\alpha} f(t)
			& = \frac{1}{g'(t)} \, \frac { d }{ dt } \left( { J }_{ T,a^+ }^{ 1-\alpha  }f (t) \right),
			\\
			D_{T,b^- }^{\alpha} f(t)
			& = \frac{-1}{g'(t)} \, \frac { d }{ dt } \left( { J }_{ T,b^- }^{ 1-\alpha  }f (t) \right).
		\end{aligned}
	\end{equation}
	for each $t \in (a,b)$.
\end{definition}

Note that if we choose
$$
g(t)=t,
\quad G( x ,\alpha ) = \Gamma(\a)\, x^{ 1-\alpha},
\quad T(t,s,\alpha )=\Gamma (\alpha )\, |t-s|^{ 1-\alpha},
$$
then $D_{T,a^+ }^{\alpha}f(t)= \, _{ }^{ RL }\!{ { D}_{a^+ }^{ \alpha  } }f(t)$
and $D_{T,b^- }^{\alpha}f(t)= \, _{ }^{ RL }\!{ { D}_{b^-}^{ \alpha  } }f(t)$.
Also, we can obtain Hadamard and others fractional derivatives
as particular cases of this generalized derivative.

\bigskip

\section{Opial-type inequality}

In 1960, Opial \cite{Opial} proved the following inequality:

\medskip

If $f \in C^1[0, h]$ satisfies $f(0) = f(h) = 0$ and $f(x) > 0$ for all $x \in (0,h)$, then
$$
\int_0^h |f (x)f'(x)| \, dx
\le \frac{h}{4} \int_0^h |f'(x)|^2 \, dx .
$$
Opial's inequality and its generalizations play a main
role in establishing the existence and uniqueness of initial and boundary value
problems for ordinary and partial differential equations \cite{2,3,4,5,6}.
For an extensive survey on these Opial-type inequalities, see \cite{2,6}.

\medskip

We need the following result in \cite[p.44]{M}
(see the original proof in \cite{Mu}).
Although the result in \cite[p.44]{M} deals with measures on $(0,\infty)$, it can be reformulated for measures on a compact interval
(see, e.g., \cite[Theorem 3.1]{APRR}).

\medskip

\noindent {\bf Muckenhoupt inequality.} {\it Let us consider $1\le p \le q < \infty$ and
measures $\mu_0,\mu_1$ on $[a,b]$ with $\mu_0(\{b\})= 0$.
Then there exists a positive constant $C$ such that
$$
\Big\|\int_{a}^x u(t)\,dt \Big\|_{L^q([a,b],\mu_0)}\le
C \,
\big\|u\big\|_{L^p([a,b],\mu_1)}
$$
for any measurable function $u$ on $[a,b]$, if and only if}
\begin{equation} \label{eq:MI2}
B := \sup_{a<x<b} \mu_0([x,b))^{1/q}
 \big\|(d\mu_1/dx)^{-1}\big\|_{L^{1/(p-1)}([a,x])}^{1/p}
< \infty ,
\end{equation}
{\it where we use the convention $0\cdot \infty=0$.
Moreover, we can choose
\begin{equation} \label{eq:MI3}
C =
\begin{cases}
B \Big( \displaystyle\frac{q}{q-1} \Big)^{(p-1)/p} q^{1/q} , & \quad \text{ if $p > 1$,}
\\
B, & \quad \text{ if $p=1$.}
\end{cases}
\end{equation}
}

\smallskip

Muckenhoupt inequality allows to improve Opial inequality in several ways:

\smallskip

$(1)$ we allow to integrate with respect to very general measures,

\smallskip

$(2)$ we do not need the hypothesis $f(b) = 0$,

\smallskip

$(3)$ we do not need the hypothesis $f > 0$ on $(a,b)$,

\smallskip

$(4)$ we substitute the hypothesis $f \in C^1[a,b]$ by the weaker one: $f$ is an absolutely continuous function on $[a,b]$.

\begin{theorem} \label{t:Opial}
Let us consider $1\le p \le q < \infty$ and
measures $\mu_0,\mu_1$ on $[a,b]$ with $\mu_0(\{b\})= 0$
and
$$
B := \sup_{a<x<b} \mu_0([x,b))^{1/q}
 \big\|(d\mu_1/dx)^{-1}\big\|_{L^{1/(p-1)}([a,x])}^{1/p}
< \infty .
$$
Then
$$
\big\| f f' \big\|_{L^1([a,b],\mu_0)}
\le B \Big( \frac{q}{q-1} \Big)^{(p-1)/p} q^{1/q} \,
\big\|f'\big\|_{L^p([a,b],\mu_1)} \big\| f' \big\|_{L^{q/(q-1)}([a,b],\mu_0)}
$$
if $p>1$, and
$$
\big\| f f' \big\|_{L^1([a,b],\mu_0)}
\le B \,
\big\|f'\big\|_{L^1([a,b],\mu_1)} \big\| f' \big\|_{L^{q/(q-1)}([a,b],\mu_0)} ,
$$
for every absolutely continuous function $f$ on $[a,b]$ with $f(a)=0$.
\end{theorem}

\begin{proof}
By Muckenhoupt inequality, the constant
$$
C =
\begin{cases}
B \Big( \displaystyle\frac{q}{q-1} \Big)^{(p-1)/p} q^{1/q} , & \quad \text{ if $p > 1$,}
\\
B, & \quad \text{ if $p=1$,}
\end{cases}
$$
satisfies
$$
\Big\|\int_{a}^x u(t)\,dt \Big\|_{L^q([a,b],\mu_0)}\le
C \,
\big\|u\big\|_{L^p([a,b],\mu_1)}
$$
for any measurable function $u$ on $[a,b]$.
For each absolutely continuous function $f$ on $[a,b]$ with $f(a)=0$, we have
that there exists $f'$ a.e. on $[a,b]$,
$f' \in L^1[a,b]$ and
$$
f(x)=\int_{a}^x f'(t)\,dt
$$
for every $x \in [a,b]$.
Consequently,
$$
\big\| f \big\|_{L^q([a,b],\mu_0)}
\le C \,
\big\|f'\big\|_{L^p([a,b],\mu_1)} .
$$
Hence, H\"older inequality gives
$$
\begin{aligned}
\big\| f f' \big\|_{L^1([a,b],\mu_0)}
& \le \big\| f \big\|_{L^q([a,b],\mu_0)} \big\| f' \big\|_{L^{q/(q-1)}([a,b],\mu_0)}
\\
& \le C \, \big\|f'\big\|_{L^p([a,b],\mu_1)} \big\| f' \big\|_{L^{q/(q-1)}([a,b],\mu_0)} .
\end{aligned}
$$
\end{proof}

\begin{remark} \label{r:Opial}
For each absolutely continuous function $f$ on $[a,b]$
the set
$$
S=\{ \, x \in [a,b] : \, \nexists \, f'(x) \,\}
$$
has zero Lebesgue measure, but it is possible to have $\mu_0(S) > 0$ and/or $\mu_1(S) > 0$.
The argument in the proof of Theorem \ref{t:Opial} gives that the inequality holds for any fixed choice of values of $f'$ on $S$.
\end{remark}

Theorem \ref{t:Opial} has the following direct consequence.

\begin{corollary} \label{c:Opial}
Let us consider $1\le p \le q < \infty$ and
a measure $\mu$ on $[a,b]$ with $\mu(\{b\})=0$
and
$$
B := \sup_{a<x<b} \mu([x,b))^{1/q}
 \big\|(d\mu/dx)^{-1}\big\|_{L^{1/(p-1)}([a,x])}^{1/p}
< \infty .
$$
Then
$$
\big\| f f' \big\|_{L^1([a,b],\mu)}
\le B \Big( \frac{q}{q-1} \Big)^{(p-1)/p} q^{1/q} \,
\big\|f'\big\|_{L^p([a,b],\mu)} \big\| f' \big\|_{L^{q/(q-1)}([a,b],\mu)}
$$
if $p>1$, and
$$
\big\| f f' \big\|_{L^1([a,b],\mu)}
\le B \,
\big\|f'\big\|_{L^1([a,b],\mu)} \big\| f' \big\|_{L^{q/(q-1)}([a,b],\mu)} ,
$$
for every absolutely continuous function $f$ on $[a,b]$ with $f(a)=0$.
\end{corollary}

Corollary \ref{c:Opial} has the following consequence.

\begin{corollary} \label{c:Opial2}
Let us consider $1\le p \le 2$ and
a measure $\mu$ on $[a,b]$ with $\mu(\{b\})=0$
and
$$
B := \sup_{a<x<b} \mu([x,b))^{(p-1)/p}
 \big\|(d\mu/dx)^{-1}\big\|_{L^{1/(p-1)}([a,x])}^{1/p}
< \infty .
$$
If $1< p \le 2$, then
$$
\big\| f f' \big\|_{L^1([a,b],\mu)}
\le B \Big( \frac{p^2}{p-1} \Big)^{(p-1)/p} \,
\big\|f'\big\|_{L^p([a,b],\mu)}^2
$$
for every absolutely continuous function $f$ on $[a,b]$ with $f(a)=0$.

Furthermore, if $\mu$ is a finite measure, then
$$
\big\| f f' \big\|_{L^1([a,b],\mu)}
\le B \,
\big\|f'\big\|_{L^1([a,b],\mu)}^2
$$
for every absolutely continuous function $f$ on $[a,b]$ such that $f(a)=0$.
\end{corollary}

\begin{proof}
Assume first that $1 < p \le 2$.
Let us consider $q \ge 2$ such that $1/p+1/q=1$, and so, $p=q/(q-1)$ and $q=p/(p-1)$.
Thus, $1 < p \le 2 \le q < \infty$ and
Corollary \ref{c:Opial} gives the result, since
$$
B \Big( \frac{q}{q-1} \Big)^{(p-1)/p} q^{1/q}
= B \, p^{(p-1)/p} \Big( \frac{p}{p-1} \Big)^{(p-1)/p}
= B \Big( \frac{p^2}{p-1} \Big)^{(p-1)/p} .
$$

Assume now that $\mu$ is a finite measure,
and fix an absolutely continuous function $f$ on $[a,b]$ such that $f(a)=0$ and $f'\in L^{p_0}([a,b],\mu)$ for some $p_0 > 1$.
We have proved that
$$
\big\| f f' \big\|_{L^1([a,b],\mu)}
\le B \Big( \frac{p^2}{p-1} \Big)^{(p-1)/p} \,
\big\|f'\big\|_{L^p([a,b],\mu)}^2
$$
for every $1 < p \le \min \{p_0,2\}$.

Let us consider $B=B(p)$ as a function of $p$.
Thus,
$$
\begin{aligned}
B(p)
& \le \mu([a,b))^{(p-1)/p}
 \big\|(d\mu/dx)^{-1}\big\|_{L^{1/(p-1)}([a,b])}^{1/p} .
\end{aligned}
$$
Since $\mu$ is a finite measure, we have
$$
\begin{aligned}
\limsup_{p \to 1^+} B(p)
& \le \lim_{p \to 1^+} \mu([a,b))^{(p-1)/p}
 \big\|(d\mu/dx)^{-1}\big\|_{L^{1/(p-1)}([a,b])}^{1/p}
\\
& = \big\|(d\mu/dx)^{-1}\big\|_{L^{\infty}([a,b])}
= B(1) .
\end{aligned}
$$
Since
$$
|f'|^p
\le |f'|^{p_0} \chi_{\{ |f'| \ge 1 \}} + \chi_{\{ |f'| < 1 \}}
\le |f'|^{p_0} + 1
\in L^1([a,b], \mu)
$$
for every $1 < p \le p_0$,
dominated convergence theorem gives
$$
\lim_{p \to 1^+} \big\|f'\big\|_{L^p([a,b],\mu)}^2
= \big\|f'\big\|_{L^1([a,b],\mu)}^2 .
$$
Finally, we have
$$
\lim_{p \to 1^+} \Big( \frac{p^2}{p-1} \Big)^{(p-1)/p}
=1,
$$
and the desired inequality holds if $f'\in L^{p_0}([a,b],\mu)$ for some $p_0 > 1$.

Let us consider now any absolutely continuous function $f$ on $[a,b]$ such that $f(a)=0$.
Define the measure $\mu^*$ on $[a,b]$ by $d\mu^* = d\mu + dx$.
Since $f$ is an absolutely continuous function on $[a,b]$,
$f' \in L^1[a,b]$.
If $f' \notin L^1([a,b],\mu)$, then the inequality is direct.
So, we can assume that $f' \in L^1([a,b],\mu)$.
Thus, there exists a sequence $\{s_n\}$ of simple functions with
$$
\lim_{n \to \infty} \big\| f' - s_n \big\|_{L^1([a,b],\mu^*)} = 0.
$$
Hence, there exists $N$ such that
$$
\| s_n \|_{L^1([a,b],\mu^*)} - \big\| f' \big\|_{L^1([a,b],\mu^*)}
\le \big\| f' - s_n \big\|_{L^1([a,b],\mu^*)}
< 1
$$
for every $n \ge N$.
Therefore,
$$
\| s_n \|_{L^1([a,b],\mu)}
\le \| s_n \|_{L^1([a,b],\mu^*)}
\le \big\| f' \big\|_{L^1([a,b],\mu^*)} +1
$$
for every $n \ge N$.

Since $\mu$ is a finite measure,
if we define $f_n(x) = \int_a^x s_n(t)\, dt$, then $f_n \in C[a,b] \subset L^p([a,b],\mu)$ for every $p \ge 1$,
and we have proved that
$$
\big\| f_n f_n' \big\|_{L^1([a,b],\mu)}
\le B \,
\big\|f_n'\big\|_{L^1([a,b],\mu)}^2 .
$$
We have for any $x \in [a,b]$
$$
\big| f(x) - f_n(x) \big|
= \Big| \int_a^x \big(f'(t)-s_n(t)\big)\, dt \Big|
\le \int_a^x \big|f'(t)-s_n(t)\big|\, dt
\le \big\| f' - s_n \big\|_{L^1([a,b],\mu^*)} .
$$
We have
$$
\begin{aligned}
& \big\| ff'-f_n f_n' \big\|_{L^1([a,b],\mu)}
= \int_a^b \big| ff'-f_n f_n' \big|\, d\mu
\\
& \qquad \le \int_a^b \big| ff'-f f_n' \big|\, d\mu + \int_a^b \big| ff_n'-f_n f_n' \big|\, d\mu
\\
& \qquad \le \| f \|_{\infty}\int_a^b \big| f'- f_n' \big|\, d\mu
+ \| f' - s_n \|_{L^1([a,b],\mu^*)} \int_a^b | s_n |\, d\mu
\\
& \qquad \le \| f \|_{\infty} \| f' - s_n \|_{L^1([a,b],\mu^*)}
+ \| f' - s_n \|_{L^1([a,b],\mu^*)} \big( \big\| f' \big\|_{L^1([a,b],\mu^*)} +1 \big)
\end{aligned}
$$
for every $n \ge N$.
Hence,
$$
\lim_{n \to \infty} \big\| ff'-f_n f_n' \big\|_{L^1([a,b],\mu)}
= 0
$$
and so,
$$
\big\| f f' \big\|_{L^1([a,b],\mu)}
\le B \,
\big\|f'\big\|_{L^1([a,b],\mu)}^2 .
$$
\end{proof}

\medskip

If we choose $\mu$ as the Lebesgue measure on $[a,b]$, then we obtain the following results.

\begin{corollary} \label{c:Opial3}
Let us consider $1\le p \le q < \infty$.
Then
$$
\big\| f f' \big\|_{L^1([a,b])}
\le \Big( \frac{b-a}{1/q+(p-1)/p} \Big)^{1/q+(p-1)/p} \Big( \frac{q(p-1)}{p(q-1)} \Big)^{(p-1)/p}
\big\|f'\big\|_{L^p([a,b])} \big\| f' \big\|_{L^{q/(q-1)}([a,b])}
$$
if $p>1$, and
$$
\big\| f f' \big\|_{L^1([a,b])}
\le ( b-a )^{1/q} \big\|f'\big\|_{L^1([a,b])} \big\| f' \big\|_{L^{q/(q-1)}([a,b])} ,
$$
for every absolutely continuous function $f$ on $[a,b]$ with $f(a)=0$.
\end{corollary}

\begin{proof}
Let us compute
$$
B = \sup_{a<x<b} (b-x)^{1/q} (x-a)^{(p-1)/p} .
$$
For each $\a>0$ and $\b \ge 0$, consider the function $u$ defined on $[a,b]$ as
$$
u(x) = (b-x)^{\a} (x-a)^{\b} .
$$
If $\b=0$, then
$$
\sup_{a<x<b} u(x)
= u(a)
= (b-a)^{\a} .
$$
Assume now that $\b > 0$.
We have for $a<x<b$
$$
\begin{aligned}
u'(x)
& = -\a(b-x)^{\a-1} (x-a)^{\b} + \b (b-x)^{\a} (x-a)^{\b-1}
= 0
\\
& \Leftrightarrow
\qquad
\b (b-x)^{\a} (x-a)^{\b-1} = \a(b-x)^{\a-1} (x-a)^{\b}
\\
& \Leftrightarrow
\qquad
\b (b-x) = \a (x-a)
\\
& \Leftrightarrow
\qquad
x = \frac{a \a + b \b}{\a+\b} \,.
\end{aligned}
$$
Since $u(a)=u(b)=0$, we have
$$
\begin{aligned}
\sup_{a<x<b} u(x)
& = \max_{a \le x \le b} u(x)
= u \Big( \frac{a \a + b \b}{\a+\b} \Big)
\\
& = \Big( \frac{\a (b-a)}{\a+\b} \Big)^{\a} \Big( \frac{\b (b-a)}{\a+\b} \Big)^{\b}
= \frac{\a^\a \b^\b}{(\a+\b)^{\a+\b}} (b-a)^{\a+\b} .
\end{aligned}
$$
Thus,
$B=( b-a )^{1/q}$ if $p=1$,
and
$$
B
= \frac{(1/q)^{1/q} ((p-1)/p)^{(p-1)/p}}{(1/q+(p-1)/p)^{1/q+(p-1)/p}} (b-a)^{1/q+(p-1)/p} ,
$$
$$
B \Big( \frac{q}{q-1} \Big)^{(p-1)/p} q^{1/q}
= \Big( \frac{b-a}{1/q+(p-1)/p} \Big)^{1/q+(p-1)/p} \Big( \frac{q(p-1)}{p(q-1)} \Big)^{(p-1)/p}
$$
if $p>1$.
Hence, Corollary \ref{c:Opial} gives the result.
\end{proof}

\begin{corollary} \label{c:Opial4}
Let us consider $1 \le p \le 2$.
Then
$$
\big\| f f' \big\|_{L^1([a,b])}
\le \Big( \frac{p(b-a)}{2(p-1)^{1/2}} \Big)^{2(p-1)/p} \,
\big\|f'\big\|_{L^p([a,b])}^2
$$
if $1<p \le 2$, and
$$
\big\| f f' \big\|_{L^1([a,b])}
\le \big\|f'\big\|_{L^1([a,b])}^2
$$
for every absolutely continuous function $f$ on $[a,b]$ such that $f(a)=0$.
\end{corollary}

\begin{proof}
Assume that $1 < p \le 2$.
It suffices to consider $q \ge 2$ such that $1/p+1/q=1$
(recall that $p=q/(q-1)$ and $q=p/(p-1)$), and to apply Corollary \ref{c:Opial3}:
$$
\begin{aligned}
& \Big( \frac{b-a}{1/q+(p-1)/p} \Big)^{1/q+(p-1)/p} \Big( \frac{q(p-1)}{p(q-1)} \Big)^{(p-1)/p}
\\
& \quad = \Big( \frac{b-a}{2(p-1)/p} \Big)^{2(p-1)/p} \Big( \frac{p(p-1)}{p} \Big)^{(p-1)/p}
\\
& \quad = \Big( \frac{p(b-a)}{2(p-1)} \Big)^{2(p-1)/p} (p-1)^{(p-1)/p}
= \Big( \frac{p(b-a)}{2(p-1)^{1/2}} \Big)^{2(p-1)/p} .
\end{aligned}
$$

Let us consider now the case $p=1$.
Since the Lebesgue measure on $[a,b]$ is finite,
Corollary \ref{c:Opial2} gives
$$
\big\| f f' \big\|_{L^1([a,b])}
\le B \,
\big\|f'\big\|_{L^1([a,b])}^2
$$
with
$$
B = \sup_{a<x<b} (b-x)^{(p-1)/p}
\| 1 \|_{L^{1/(p-1)}([a,x])}^{1/p}
= \sup_{a<x<b} (b-x)^{0}
\| 1 \|_{L^{\infty}([a,x])}
 = 1.
$$
\end{proof}

\begin{remark}
Note that in the second inequality in Corollary \ref{c:Opial4}:
$$
\big\| f f' \big\|_{L^1([a,b])}
\le \big\|f'\big\|_{L^1([a,b])}^2,
$$
the constant $1$ multiplying $\|f' \|_{L^1([a,b])}^2$
does not depend on the length of the interval $[a,b]$.
\end{remark}

Corollaries \ref{c:Opial} and \ref{c:Opial2}
have, respectively, the following direct consequences for general fractional integrals of Riemann-Liouville type.

\begin{proposition} \label{p:OpialRL}
Let us consider $1\le p \le q < \infty$. If
$$
B := \sup_{a<x<b} \Big( \int _{ x }^{ b } \frac { 1 }{ T(b,s,\alpha ) } \, ds \Big)^{1/q}
 \Big( \int _{ a }^{ x } T(b,s,\alpha )^{1/(p-1)} \, ds \Big)^{(p-1)/p}
< \infty ,
$$
then
$$
\int_{a}^b \frac { | f(s) f'(s) | }{ T\big(b,s,\alpha\big) } \, ds
\le B \Big( \frac{q}{q-1} \Big)^{(p-1)/p} q^{1/q} \,
\Big( \int_{a}^b \frac { | f'(s) |^p }{ T\big(b,s,\alpha\big) } \, ds \Big)^{1/p}
\Big( \int_{a}^b \frac { | f'(s) |^{q/(q-1)} }{ T\big(b,s,\alpha\big) } \, ds \Big)^{(q-1)/q}
$$
if $p>1$, and
$$
\int_{a}^b \frac { | f(s) f'(s) | }{ T\big(b,s,\alpha\big) } \, ds
\le B \int_{a}^b \frac { | f'(s) | }{ T\big(b,s,\alpha\big) } \, ds \,
\Big( \int_{a}^b \frac { | f'(s) |^{q/(q-1)} }{ T\big(b,s,\alpha\big) } \, ds \Big)^{(q-1)/q}
$$
for every absolutely continuous function $f$ on $[a,b]$ with $f(a)=0$.
\end{proposition}

\begin{proposition} \label{p:OpialRL2}
Let us consider $1 \le p \le 2$ and assume that
$$
B := \sup_{a<x<b} \Big( \int _{ x }^{ b } \frac { 1 }{ T(b,s,\alpha ) } \, ds \Big)^{(p-1)/p}
 \Big( \int _{ a }^{ x } T(b,s,\alpha )^{1/(p-1)} \, ds \Big)^{(p-1)/p}
< \infty .
$$
If $1 < p \le 2$, then
$$
\int_{a}^b \frac { | f(s) f'(s) | }{ T\big(b,s,\alpha\big) } \, ds
\le B \Big( \frac{p^2}{p-1} \Big)^{(p-1)/p}
\Big( \int_{a}^b \frac { | f'(s) |^p }{ T\big(b,s,\alpha\big) } \, ds \Big)^{2/p}
$$
for every absolutely continuous function $f$ on $[a,b]$ with $f(a)=0$.

Furthermore, if
$$
\int_{a}^b \frac { ds }{ T\big(b,s,\alpha\big) } < \infty,
$$
then
$$
\int_{a}^b \frac { | f(s) f'(s) | }{ T\big(b,s,\alpha\big) } \, ds
\le B
\Big( \int_{a}^b \frac { | f'(s) | }{ T\big(b,s,\alpha\big) } \, ds \Big)^{2}
$$
for every absolutely continuous function $f$ on $[a,b]$ with $f(a)=0$.
\end{proposition}

\section*{Acknowledgments}

The research all the authors is
	supported by a grant from Agencia Estatal de Investigaci\'on (PID2019-106433GB-I00 / AEI / 10.13039/501100011033), Spain.
The research of Jos\'e M. Rodr\'iguez is supported by the Madrid Government (Comunidad de Madrid-Spain) under the Multiannual Agreement with UC3M in the line of Excellence of University Professors (EPUC3M23), and in the context of the V PRICIT (Regional Programme of Research and Technological Innovation).

\end{document}